\documentclass[preprint,12pt]{elsarticle}
\usepackage{etex} 

\usepackage[english]{babel}
\usepackage{amsmath}
\usepackage[standard,thmmarks,amsmath]{ntheorem}
\usepackage{amssymb}
\usepackage{graphicx}
\usepackage[normalem]{ulem}
\usepackage{algorithm}
\usepackage{algpseudocode}
\usepackage{enumerate}
\usepackage{subfigure}
\usepackage{lineno}
\usepackage{hyperref}

\journal{Applied Mathematics and Computation}

\def\le{\leqslant}

\def\ge{\geqslant}

\usepackage{amsmath} 
 
 \DeclareMathOperator*{\argmax}{arg\,max}

\begin{document}

\begin{frontmatter}



\title{Maximizing Wiener Index for Trees with Given Vertex Weight and Degree Sequences\tnoteref{mytitlenote}}
\tnotetext[mytitlenote]{This research is supported by the grant of Russian Foundation for Basic Research, project No 16--37--60102
mod\_a\_dk.}

\author{Mikhail Goubko}

\address{65 Profsoyuznaya str., 117997, Moscow, Russia\corref{mycorrespondingauthor}}

\begin{abstract}
The Wiener index is maximized over the set of trees with the given vertex weight and
degree sequences. This model covers the traditional ``unweighed'' Wiener index, the
terminal Wiener index, and the vertex distance index. It is shown that there exists
an optimal caterpillar. If weights of internal vertices increase in their degrees,
then an optimal caterpillar exists with weights of internal vertices on its backbone
monotonously increasing from some central point to the ends of the backbone, and the
same is true for pendent vertices. A tight upper bound of the Wiener index value is
proposed and an efficient greedy heuristics is developed that approximates well the
optimal index value. Finally, a branch and bound algorithm is built and tested for
the exact solution of this NP-complete problem.
\end{abstract}


\end{frontmatter}


\section{Nomenclature}\label{section_nomenclature}
This section introduces the basic graph-theoretic notation. The \textit{vertex set} and the \textit{edge set}
of a simple connected undirected graph $G$ are denoted with $V(G)$ and $E(G)$
respectively, and the \textit{degree} (i.e., the number of incident edges) of vertex $v\in V(G)$ in graph $G$ is denoted with $d_G(v)$. Let $W(G)$ be the set of \emph{pendent vertices} (those having degree one)
of graph $G$, and let $M(G):=V(G)\backslash W(G)$ be the set of its \emph{internal
vertices}. Connected graph $T$ with $|E(T)|=|V(T)|-1$ is called a \emph{tree}. Let
$\mathcal{T}$ denote the set of all trees.
\begin{definition}
A tree is called a \textit{path} if it has exactly two pendent
vertices.
\end{definition}
\begin{definition}
A tree is a \textit{caterpillar} if removing pendent
vertices and their incident edges makes a path (called the
\textit{backbone} of this caterpillar).
\end{definition}
\begin{definition}
A tree is called a \textit{star} if it has at most one internal vertex.
\end{definition}
\begin{definition}
In a \textit{starlike tree} the degree of at most one vertex exceeds 2.
\end{definition}
\begin{definition}
The \textit{centroid} is a midpoint of the longest path in the tree.
\end{definition}

Graph $G$ is called \textit{vertex-weighted} if non-negative weight is assigned to
each its vertex. The weight of vertex $v\in V(G)$ in graph $G$ is denoted as
$\mu_G(v)$. Let $\mathcal{WT}$ stand for the set of all vertex-weighted trees.

Let us consider monotone decreasing natural sequence $d_i, i=1,...,n$, and non-negative sequence $\mu_i,
i=1,...,n$, of the same length and introduce corresponding column vectors
$\mathbf{d}=(d_1, ..., d_n)^\mathrm{T}$, $\mathbf{w}=(\mu_1,...,\mu_n)^\mathrm{T}$. 

\begin{definition}
Tree $T$ \textit{has degree sequence} $\mathbf{d}$ if its vertices can be indexed
from $v_1$ to $v_n$ such that $d_T(v_i)=d_i$, $i=1,..., n$. Let
$\mathcal{T}(\mathbf{d})$ be the set of trees with degree sequence $\mathbf{d}$.
\end{definition}
It is known that $\mathcal{T}(\mathbf{d})$ is not empty, if and only if
\begin{equation}\label{eq_tree_degrees}
d_1+...+d_n=2(n-1).
\end{equation}

\begin{definition}
Vertex-weighted tree $T\in \mathcal{WT}$ \textit{has vertex degree sequence
$\mathbf{d}$ and vertex weight sequence $\mathbf{w}$} if its vertices can be indexed
from $v_1$ to $v_n$ such that $d_T(v_i)=d_i, \mu_T(v_i)=\mu_i$ for all $i=1,...,n$.
Let $\mathcal{WT}(\mathbf{w}, \mathbf{d})$ be the set of all such trees.
\end{definition}
Without loss of generality assume that if $d_i=d_j$ and $i<j$, then $\mu_i \ge \mu_j$.
\begin{definition}
Weight sequence $\mathbf{w}$ is \emph{monotone} in degree sequence $\mathbf{d}$, if
from $d_i, d_j\ge 2$, $i<j$, it follows that $\mu_i \ge \mu_j$.  
\end{definition}

For any pair of vertices $u,v\in V(G)$ of connected graph $G$ let $d_G(u, v)$ be the
distance (the number of edges in the shortest path) between vertices $u$ and $v$ in graph
$G$. Then \textit{the Wiener index} of graph $G$ is defined as
\begin{equation}\label{eq_WI_def}
WI(G):=\frac{1}{2}\sum_{u,v\in V(G)}d_G(u,v).
\end{equation}

\section{Introduction}

Graph invariants (also known as topological indices) play an important role in
algebraic graph theory providing numeric measures for various structural properties
of graphs. The Wiener index (\ref{eq_WI_def}) is probably the most renowned graph invariant. It measures ``compactness'' of
a connected graph \cite{dobrynin2001wiener}; for instance, a \textit{star} has the minimum value of the Wiener index among all trees of the given order, while a \textit{path} has the
maximum value of WI. The most ``compact'' (i.e., the one minimizing the Wiener
index) tree with the given vertex degree sequence is a ``greedy'' balanced tree, in
which all distances from leaves to the \textit{centroid} differ by at most unity while vertex degrees do not decrease towards the centroid \cite{wang2008extremal,zhang2008wiener}.

The \emph{Wiener index for graphs with weighted vertices} was proposed in
\cite{klavvzar1997wiener}. It can be defined as
\begin{equation}\label{def_VWWI}
VWWI(G):=\frac{1}{2}\sum_{u,v\in V(G)}\mu_G(u) \mu_G(v) d_G (u,v),
\end{equation}
where $d_G(u, v)$ is the distance between vertices $u$ and $v$ in graph $G$, while  $\mu_G(u)$ and  $\mu_G(v)$
are, respectively, real weights of graph vertices $u$ and $v$.

VWWI is used to foster calculation of the Wiener index \cite{kelenc2015edge}, to
predict boiling and melting points of various compounds
\cite{gao2015vertex,goubkomiloserdov2016wiener}. In particular, in
\cite{goubkomiloserdov2016wiener} the search of an alcohol isomer with the minimum
normal boiling point was reduced to the minimization of VWWI over the set of trees with
the given vertex weight and degree sequences.

It is shown in \cite{goubko2016wiener} that if weights of internal vertices do not
decrease in their degrees, then the most ``compact'' (the one minimizing VWWI) tree
with the given vertex weight and degree sequences is the, so called, \textit{generalized
Huffman tree}. It is efficiently constructed by joining sequentially sub-graphs of
the minimum weight.

The problem of the Wiener index maximization appeared a bit more complex. It is
known that an extremal tree is some caterpillar \cite{shi1993average} (i.e., a tree
that makes a path, called a \textit{backbone}, after deletion of all its pendent
vertices); vertex degrees first do not increase and then do not decrease while one
moves from one to the other end of the backbone \cite{wang2008extremal}. An efficient
dynamic programming algorithm assigns internal vertices to positions on the backbone
of an optimal caterpillar \cite{ccela2011wiener}.

In this article the problem of the maximum Wiener index over the set of trees
with the given vertex weight and degree sequences is solved for the case of internal
vertex weights being monotone in degrees. This problem appears NP-complete (the
classic partition problem reduces to its special case). It is shown that, similarly
to the partition problem, complexity of the maximization problem for WI and VWWI is
a result of asymmetry of the vertex set. If for each distinct combination of the weight
and the degree the number of vertices having this weight and degree is even (although
one internal and/or one pendent vertex with minimum weight may be unmatched), then
VWWI is maximized by a symmetric caterpillar, in which vertices are placed
mirror-like with respect to its center in the order of increasing weights. For the
general case an analytical upper bound is proposed, the greedy heuristic algorithm
and the economic branch and bound scheme are constructed, and their performance is
evaluated for random weight and degree sequences.

\section{Literature Review}\label{section_literature}

Since its appearance in 1947 \cite{wiener1947structural} the Wiener index remains one
of the most discussed graph invariants. On the one hand, its mathematical properties
have been comprehensively studied (see surveys in
\cite{dobrynin2001wiener,knor2016mathematical}). On the other hand, its relation is
established to many physical and chemical properties of compounds of different
classes (\cite{diudea1998wiener,gutman1993some,rouvray1987modeling,lang2016novel} and many
others). Many papers that appeared in recent decades investigate extremal graphs that
deliver the minimum or the maximum of the Wiener index over various sets of graphs
\cite{fischermann2002wiener,lin2016segment,lin2014extremal,wang2008extremal,zhang2008wiener,krnc2016wiener,ramane2016note}
along with its lower and upper bounds. In particular, a relation is established
between the Wiener index and the Randi{\'c} index \cite{shi2016chemical}, the spectrum of the Laplacian matrix
\cite{merris1989edge,merris1990distance,mohar1991eigenvalues,mohar1991laplacian}, 
and the distance matrix (\cite{merris1990distance,indulal2009sharp} and others) of
the graph.

Graphs with prescribed vertex degrees have been studied for more than 50 years
\cite{gallai1960graphs,burr1988extremal,shi1993average,biyikoglu2008graphs}. In
particular, the trees with the given vertex degree sequence appear in extremal problems
for linear combinations of distance-based topological indices (e.g., the Wiener
index) and degree-based indices \cite{gutman2013degree} (e.g., the first Zagreb
index \cite{gutman2004first,gutman2015graphs} and its extensions
\cite{das2015zagreb,ghorbani2012new,eliasi2012multiplicative,fath2011old,xu2014maximizing,furtula2015forgotten,su2015graphs,xu2012unified}).

Wang \cite{wang2008extremal} and Zhang et al. \cite{zhang2008wiener} have shown
independently that the minimizer of the Wiener index over the set of trees with the
given sequence of vertex degrees is the, so-called, \emph{greedy tree}
\cite{wang2008extremal}. It is efficiently built in the top-down manner by adding
vertices from the highest to the lowest degree to the seed (a vertex of maximum
degree) to keep the tree as balanced as possible. The proof in
\cite{zhang2008wiener} employs the majorization theory, which has shown to be
useful in the broad variety of topological index optimization problems
\cite{zhang2015extremal,liu2015recent}.

The problem of the maximum Wiener index for trees with the given degree sequences
appeared more complex. Schmuck et al. \cite{schmuck2012greedy} have shown that an
optimal tree is a caterpillar with vertex degrees non-increasing from the ends of
the caterpillar towards its central part. An efficient algorithm was suggested in
\cite{ccela2011wiener} that optimally assigns internal vertices to the backbone of a
caterpillar.

The Wiener index for vertex-weighted graphs (denoted as $VWWI$, the
``vertex-weighted Wiener index'') is defined by expression (\ref{def_VWWI})
\cite{klavvzar1997wiener}. The \textit{terminal Wiener index}
\cite{gutman2009terminal,gutman2010survey} is obtained as its special case by
assigning zero weights to internal vertices and unit weights to pendent vertices.
The \textit{transmission} or the \textit{vertex distance index}
\cite{mohar1988compute} is the sum of distances in graph $G$ from all
vertices to the given vertex (e.g., $v_1\in V(G)$). It can be seen as a limit of
VWWI for $\mu_1=a$, $\mu_i=\frac{1}{a}, i=2,...,n$, for $a\Rightarrow \infty$
\cite{goubkomiloserdov2016wiener}). Also, VWWI clarifies the relation between the
Wiener index and spectral properties of the \textit{graph distance matrix}. If vector
$\mathbf{w}$ of vertex weights lies on the unit sphere, $VWWI(G)$ is equal to the half of
the Raileigh quotient \cite{spielman2012} for graph distance matrix $D(G):=(d_G(u,v))_{u,v\in V(G)}$ and
vector $\mathbf{w}$: $VWWI(G)=\frac{1}{2}\mathbf{w}^\mathrm{T} D(G) \mathbf{w}$. Unfortunately, VWWI is still understudied at the moment.

The Wiener index for vertex-weighted trees and the set of trees with the given vertex
weight and degree sequences play the central role in this article. It is known that
the Wiener index and VWWI (and the terminal Wiener index as its special case) have
similar properties. In particular, the majorization
theory was used in \cite{goubko2016wiener} to minimize VWWI over trees with the given vertex weight and degree
sequences. It was shown that if weights are monotone in degrees, then VWWI is
minimized with the, so-called, \textit{generalized Huffman tree}, which can be
efficiently (with complexity $O(|V| \ln |V|))$ built with some extension of the
famous Huffman algorithm for the optimal prefix code \cite{huffman1952method}). The
greedy tree mentioned above is a special case of the generalized Huffman tree for
equal vertex weights.

In this article it is shown that, similarly to the ``unweighed'' Wiener index, a
tree with the given vertex weight and degree sequences, which maximizes VWWI, is a
caterpillar. If, in addition, vertex weights are degree-monotone, then weights of
internal vertices (and also their degrees, by monotonicity of weights) first do not
increase and then do not decrease while one moves along the backbone of an optimal
caterpillar. Similarly, the weights of pendent vertices being adjacent to vertices
of the backbone first do not increase and then do not decrease while one moves along
the backbone of an optimal caterpillar. Distinct to maximization of the Wiener index
for general trees, maximization of VWWI is generally NP-complete; the
\textit{Wiener-type quadratic assignment problem} (QAP) reduces to it, while it is
known that the classic partition problem reduces to the Wiener-type QAP
\cite{ccela2011wiener}.

A quasi-polynomial dynamic programming algorithm for the Wiener-type QAP was
proposed in \cite{ccela2011wiener}. In this article the continuous relaxation of QAP
is used to propose the upper bound estimate of VWWI. A
branch and bound scheme is also constructed, which is applicable to degree-monotone weights.

\section{Caterpillar is an Optimal Tree}\label{section_caterpillar}

Let us consider the problem of the Wiener index maximization over a set of trees
with the given vertex weight and degree sequences. Hereinafter, the tree delivering the maximum
of VWWI is called an \textit{optimal} tree.

In this section it is shown that an optimal tree with positive vertex weights is
essentially a caterpillar. For the ``classical'' Wiener index this result has been
proven in \cite{shi1993average} but we follow the line of the proof from
\cite{schmuck2012greedy} instead. Let us introduce an auxiliary result.

\begin{figure}
\subfigure[Starlike
tree]{\includegraphics[width=0.5\textwidth]{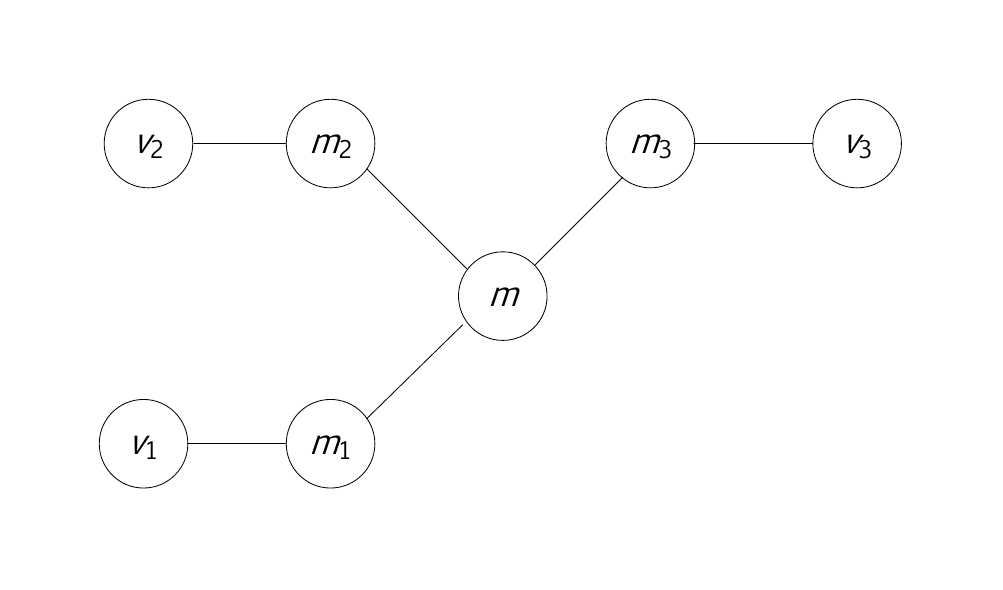}\label{fig_starlike_i}}
\subfigure[Transformed
tree]{\includegraphics[width=0.5\textwidth]{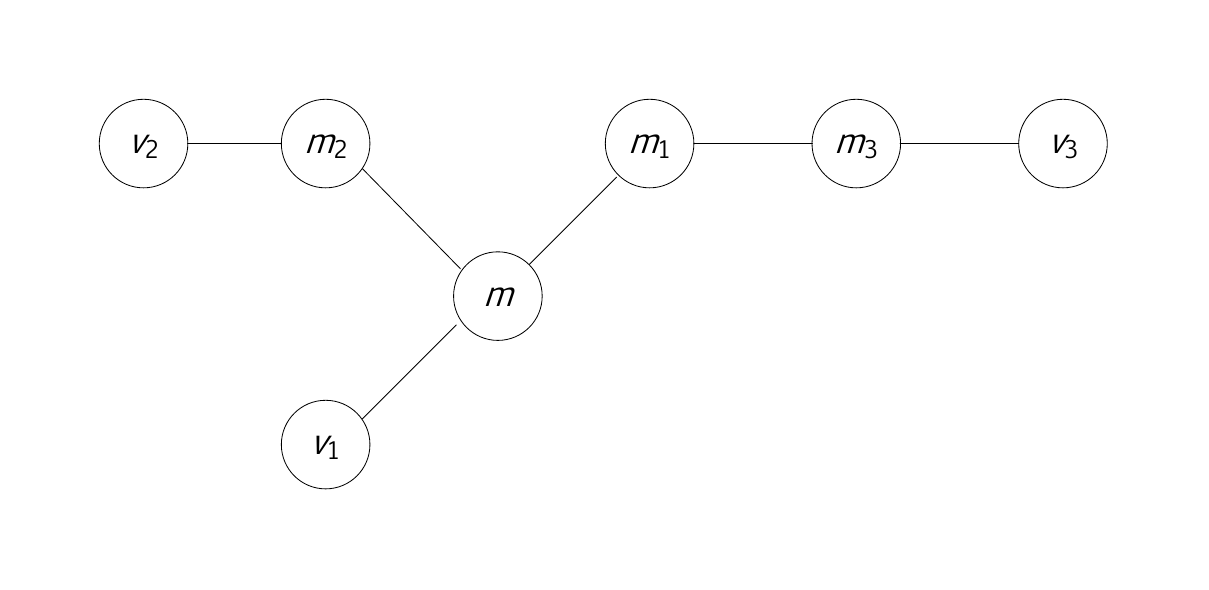}\label{fig_starlike_m}}
\caption{Transformation of stralike tree}\label{fig_starlike}
\end{figure}

\begin{lemma}\label{lemma_starlike}
Let us consider a three-leaf balanced starlike tree $T\in \mathcal{WT}$ of order 7
with \textbf{positive} vertex weights and vertices labeled as in Figure
\ref{fig_starlike_i}, and tree $T'$ obtained from $T$ by replacing edges $mm_3$ and
$m_1v_1$ with edges $mv_1$ and $m_1m_3$ as in Figure \ref{fig_starlike_m}. Let us
denote for short $\mu:=\mu_T(m), \mu_i':=\mu_T(m_i), \mu_i'':=\mu_T(v_i), i=1,2,3$.
If
\begin{equation}\label{ineq_weights}
\mu_1'+\mu_1'' \le \mu_2'+\mu_2'' \le\mu_3'+\mu_3'',
\end{equation}
then $VWWI(T)< VWWI(T')$.
\begin{proof}\sloppy
From expression (\ref{def_VWWI}), $VWWI(T)=\frac{1}{2}\mathbf{m}^\textrm{T} D(T)
\mathbf{m}$, where $\mathbf{m}:=(\mu,\mu_1',\mu_2',\mu_3',\mu_1'',\mu_2'',\mu_3'')$
is the vector of vertex weights and $D(T)$ is the distance matrix of tree $T$.
Hence, $VWWI(T')-VWWI(T)=\frac{1}{2}\mathbf{m}^\textrm{T} (D(T')-D(T)) \mathbf{m}$.

Direct calculation of distance matrices gives
\begin{equation}\label{eq_VWWIdiff}
VWWI(T')-VWWI(T)=(\mu+\mu_2'+\mu_2''-\mu_1')(\mu_3'+\mu_3''-\mu_1'').
\end{equation}
Vertex weights are positive, so
$\mu+\mu_2'+\mu_2''-\mu_1'>\mu+\mu_2'+\mu_2''-\mu_1'-\mu_1''\ge 0$, and
$\mu_3'+\mu_3''-\mu_1''> \mu_3'+\mu_3''-\mu_1'-\mu_1''\ge 0$. Therefore, $VWWI(T)<
VWWI(T')$.
\end{proof}
\end{lemma}

\begin{lemma}\label{lemma_caterpillar}
Given sequence of \textbf{positive} vertex weights $\mathbf{w}$ and degree sequence
$\mathbf{d}$ of the same length, any optimal tree $T^*\in
\mathcal{WT}(\mathbf{w},\mathbf{d})$ is a caterpillar.
\begin{proof}
All the trees with at most three internal vertices are caterpillars, so, without the loss of
generality assume that there are at least four elements with $d_i>1$. Assume, by
contradiction, that $T^*$ is not a caterpillar. Then an internal vertex $m\in
M(T^*)$ exists being adjacent to three other internal vertices (let us denote them
with $m_1, m_2$ и $m_3$). Since vertex $m_i$ is internal, it has at least one more
adjacent vertex distinct from $m$. Denote this adjacent vertex with $v_i$,
$i=1,2,3$.

Deleting any vertex $v\in V(T^*)$ and its incident edges in tree $T^*$, one obtains
$d_{T^*}(v)$ connected components called \textit{induced subtrees}. Let $\mu$ stand
for the sum of weights of vertex $m$ and the vertices of all subtrees induced by
deletion of $m$ from $T^*$ and containing neither of the vertices $m_1, m_2, m_3$.
Similarly, let $\mu_i'$, $i=1,2,3$, stand for the sum of weights of vertex $m_i$ and
vertices of all subtrees induced by deleting $m_i$ from $T^*$ and containing neither
$m$ nor $v_i$. Finally, let $\mu_i''$, $i=1,2,3$, stand for the sum of vertex
weights in a subtree induced by deletion of vertex $m_i$ and containing vertex
$v_i$.

Without loss of generality assume that $\mu_1'+\mu_1'' \le \mu_2'+\mu_2''
\le\mu_3'+\mu_3''$. Let us consider a tree $T^{**}\in
\mathcal{WT}(\mathbf{w},\mathbf{d})$ obtained from $T^*$ by replacing edges $mm_3$
and $m_1v_1$ with edges $mv_1$ and $m_1m_3$.

It is clear that $VWWI(T^{**})-VWWI(T^*)=VWWI(T')-VWWI(T)$, where $T$ and $T'$ are
the trees from the statement of Lemma \ref{lemma_starlike}. Since all vertex weights
are assumed positive, from Lemma \ref{lemma_starlike} it follows that
$VWWI(T')>VWWI(T)$ and $T^*$ cannot be optimal. The obtained contradiction completes
the proof.
\end{proof}
\end{lemma}

Lemma \ref{lemma_caterpillar} does not cover an important case of the terminal
Wiener index. If zero weights are allowed, an optimal caterpillar still exists,
though, similarly to the results by \cite{schmuck2012greedy}, not all optimal trees
have to be caterpillars.

\begin{corollary}\label{corollary_caterpillarzero} If zero weights are allowed
in Lemma \ref{lemma_caterpillar}, an optimal caterpillar exists.
\begin{proof}
Introduce the notation 
$$VWWI^*:=\max\left\{VWWI(T): T\in
\mathcal{WT}(\mathbf{w},\mathbf{d})\right\},$$ 
$$VWWI':=\max\left\{VWWI(T): T\in
\mathcal{WT}(\mathbf{w},\mathbf{d}), VWWI(T)<VWWI^*\right\}, $$
$$\delta :=
VWWI^*-VWWI'.$$
Since the set $\mathcal{WT}(\mathbf{w},\mathbf{d})$ is finite,
$\delta>0$. 

\sloppy Select any positive $\alpha< \sqrt{\frac{2\delta}{n^3}}$ and consider weight
sequence $\mathbf{w}^\alpha=(\mu^\alpha_1, ...,
\mu^\alpha_n)=(\mu_1+\alpha,...,\mu_n+\alpha)$. Let $T^\alpha$ be an arbitrary
optimal tree over $\mathcal{WT}(\mathbf{w}^\alpha,\mathbf{d})$. From Lemma
\ref{lemma_caterpillar}, $T^\alpha$ is a caterpillar. Introduce caterpillar
$T^\alpha_0\in \mathcal{WT}(\mathbf{w},\mathbf{d})$ obtained from $T^\alpha$ by
decreasing all vertex weights by $\alpha$. It is clear that
\begin{multline}
VWWI^*\le VWWI(T^\alpha)=\\
VWWI(T^\alpha_0)+\frac{\alpha^2}{2} \sum_{u,v\in V(T^\alpha)}d_{T^\alpha}(u,v)\le VWWI(T^\alpha_0)+\frac{\alpha^2 n^3}{2}.\end{multline}
Therefore, $VWWI(T^\alpha_0)\ge VWWI^*-\frac{\alpha^2 n^3}{2}>VWWI^*-\delta=VWWI'$.
Consequently, $VWWI(T^\alpha_0)= VWWI^*$, and caterpillar $T^\alpha_0$ is optimal.
\end{proof}
\end{corollary}

\section{Structure of Optimal Caterpillar}\label{section_structure}

\begin{definition}
Numeric sequence $a_1, ..., a_q$ is called \textit{V-shaped} if there exists such
$\underline{k}\in \{1,...,q\}$, that $a_k\ge a_{k+1}$ for $k< \underline{k}$ and
$a_k\le a_{k+1}$ for $k \ge \underline{k}$.
\end{definition}

\begin{definition}
Let us consider a (vertex-weighted) caterpillar $T\in \mathcal{WT}$ with backbone
$v_1, ..., v_q$. Vertex $v\in V(T)$ is \textit{associated} with backbone position
$k\in \{1,...,q\}$, if either $v=v_k$ or $v\in W(T), vv_k\in E(T)$. Let $A_T(k)$
stand for the set of vertices associated with backbone position $k$ of caterpillar
$T$ and their total weight is denoted with
\begin{equation}\label{eq_assweight}
w_T(k):=\sum_{v\in A_T(k)}\mu_T(v).
\end{equation}
The value $p_T(k):=\sum_{l=1}^q w_T(l)|k-l|$, $k=1,...,q$, is called the
\textit{price} of $k$-th position on the backbone of caterpillar $T$.
\end{definition}

It is shown in \cite{ccela2011wiener} that if caterpillar $T\in
\mathcal{T}(\mathbf{d})$ with backbone $v_1, ..., v_q$ maximizes $WI$, then sequence
$d_T(v_k)$, $k=1,...,q$, is V-shaped. This result is extended below by proving that the sequences
of pendent and internal vertex weights and the sequence of internal vertex degrees
are V-shaped in \textbf{some} optimal caterpillar, and the same is true for
\textbf{any} optimal caterpillar with positive vertex weights.

For any pair of pendent vertices $u_k\in A_T(k), v_l\in A_T(l)$ of vertex-weighted
caterpillar $T$ $d_T(u,v)=2+|k-l|$; for pendent vertex $u\in A_T(k)$ and internal
vertex $v\in A_T(l)$ $d_T(u,v)=1+|k-l|$. Therefore, the value of the Wiener index
can be written as
\begin{equation}\label{eq_VWWIass}
VWWI(T)=\sum_{k=1}^q w_T(k)p_T(k)+\sum_{v\in V(T)}\mu_T(v)\sum_{v\in
M(T)}\mu_T(v)-\sum_{v\in W(T)}\mu_T(v)^2.
\end{equation}
Note that if a pair of caterpillars share the same weight and degree sequences, they
differ only in the first term in (\ref{eq_VWWIass}).

\begin{lemma}\label{lemma_Vshapedpendent}
Let $\mathbf{w}$ be a \textbf{positive} weight sequence. Consider an optimal
caterpillar $T\in \mathcal{WT}(\mathbf{w}, \mathbf{d})$ and let $v_1, ...,v_q$ be
its backbone. If $u_k\in A_T(k)$ is an arbitrary pendent vertex associated with
$k$-th backbone position, then the sequence of pendent vertex weights $\mu_T(u_k)$,
$k=1,...,q$, is V-shaped.
\begin{proof}
Assume, by contradiction, that the sequence of pendent vertex weights is not
V-shaped. Then $\mu_T(u_{k-1})<\mu_T(u_k)>\mu_T(u_{k+1})$ for some
$k\in\{2,...,q-1\}$. Let us consider caterpillar $T'\in \mathcal{WT}(\mathbf{w},
\mathbf{d})$ obtained from $T$ by replacing edges $u_kv_k$ and $u_{k-1}v_{k-1}$ with
edges $u_kv_{k-1}$ and $u_{k-1}v_k$ and also caterpillar $T''\in
\mathcal{WT}(\mathbf{w}, \mathbf{d})$ obtained from $T$ by replacing edges $u_kv_k$
and $u_{k+1}v_{k+1}$ with edges $u_kv_{k+1}$ and $u_{k+1}v_k$. From Equation
(\ref{eq_VWWIass}) it follows that
$$VWWI(T')-VWWI(T)=\frac{1}{2}\left(\mu_T(u_k)-\mu_T(u_{k-1})\right)\left(p_{T+T'}(k-1)-p_{T+T'}(k)\right),$$
$$VWWI(T'')-VWWI(T)=\frac{1}{2}\left(\mu_T(u_k)-\mu_T(u_{k+1})\right)\left(p_{T+T''}(k+1)-p_{T+T''}(k)\right),$$
where $p_{T+T'}(k)$ is shorthand for $p_T(k)+p_{T'}(k)$. Since $T$ is optimal,
$VWWI(T')\le VWWI(T) \ge VWWI(T'')$. By assumption,
$\mu_T(u_{k-1})<\mu_T(u_k)>\mu_T(u_{k+1})$, so
$$p_{T+T'}(k-1)\le p_{T+T'}(k),
p_{T+T''}(k+1)\le p_{T+T''}(k)$$ and, therefore,
\begin{equation}\label{eq_concave}
p_{T}(k-1)+p_{T'}(k-1)+p_{T}(k+1)+p_{T''}(k+1)\le 2p_T(k) + p_{T'}(k)+p_{T''}(k).
\end{equation}

On the other hand, by the definition of position prices,
$$p_{T'}(k)-p_T(k)=p_T(k-1)-p_{T'}(k-1)= \mu_T(u_k)-\mu_T(u_{k-1}), $$
$$p_{T''}(k)-p_T(k)=p_T(k+1)-p_{T''}(k+1) = \mu_T(u_k)-\mu_T(u_{k+1}).$$
Substituting into Expression (\ref{eq_concave}) and dividing by two, we have
$$p_{T}(k-1)+p_{T}(k+1)\le 2p_T(k) +\mu_T(u_k)-\mu_T(u_{k-1})+\mu_T(u_k)-\mu_T(u_{k+1}).$$
Taking into account that $p_T(k-1)+p_T(k+1)-2p_T(k)=2w_T(k)\ge 2\mu_T(u_k)$, we
finally obtain $\mu_T(u_{k-1})+\mu_T(u_{k+1})\le 0$, which is impossible for
positive vertex weights. The obtained contradiction completes the proof.
\end{proof}
\end{lemma}

\begin{lemma}\label{lemma_Vshapedinternal}
Assume \textbf{positive} weight sequence $\mathbf{w}$ is monotone in degrees
$\mathbf{d}$. Then in optimal caterpillar $T\in \mathcal{WT}(\mathbf{w},
\mathbf{d})$ with backbone $v_1, ...,v_q$, the sequence of internal vertex weights
$\mu_T(v_k)$ and the sequence of internal vertex degrees $d_T(v_k)$, $k=1,...,q$,
are V-shaped.
\begin{proof}
Assume that the sequence of internal vertex weights is not V-shaped. Then
$\mu_T(v_{k-1})<\mu_T(v_k)>\mu_T(v_{k+1})$ for some $k\in\{2, ..., q-1\}$. Since
weights are degree-monotone, we also conclude that $d_T(v_{k-1})\le d_T(v_k)\ge
d_T(v_{k+1})$. Denote $d' := d_T(v_k) - d_T(v_{k-1})\ge 0$ and $d'' := d_T(v_k) -
d_T(v_{k+1})\ge 0$.

Let us select $d'$ arbitrary pendent vertices adjacent to vertex $v_k$ and denote
their total weight with $\mu'\ge 0$. Consider caterpillar $T'$ obtained from $T$ by
reconnecting all the other incident edges of vertex $v_k$ to vertex $v_{k-1}$ and
reconnecting all the incident edges of vertex $v_{k-1}$ to vertex $v_k$. It is clear
that $T'\in \mathcal{WT}(\mathbf{w},\mathbf{d})$ and
$$VWWI(T')-VWWI(T)=\frac{1}{2}(\mu'+\mu_T(v_k)-\mu_T(v_{k-1}))\left(p_{T+T'}(k-1)-p_{T+T'}(k)\right).$$

Caterpillar $T$ is optimal, so $VWWI(T')\le VWWI(T)$ and since, by assumption,
$\mu_T(v_k)>\mu_T(v_{k-1})$, $\mu'\ge 0$, we conclude that $p_{T+T'}(k-1)\le
p_{T+T'}(k)$. In a similar manner it is shown that $p_{T+T''}(k+1)\le p_{T+T''}(k)$
and, therefore,
\begin{equation}\label{eq_case2}
p_{T}(k-1)+p_{T'}(k-1)+p_{T}(k+1)+p_{T''}(k+1)\le 2p_T(k) + p_{T'}(k)+p_{T''}(k).
\end{equation}

Similarly to Lemma \ref{lemma_Vshapedpendent}, we write
$$p_{T'}(k)-p_T(k)=p_T(k-1)-p_{T'}(k-1)= \mu_T(v_k)-\mu_T(v_{k-1})+\mu', $$
$$p_{T''}(k)-p_T(k)=p_T(k+1)-p_{T''}(k+1) = \mu_T(v_k)-\mu_T(v_{k+1})+\mu',$$
$$p_T(k-1)+p_T(k+1)-2p_T(k)=2w_T(k)\ge 2(\mu_T(v_k)+\mu')>0,$$
so inequality (\ref{eq_case2}) reduces to $\mu_T(v_{k-1})+\mu_T(v_{k+1})\le 0$,
which is impossible for positive weights. The obtained contradiction proves that the
sequence of internal vertex weights is V-shaped. The same argument can be used to
show that the sequence of internal vertex degrees is, indeed, V-shaped.
\end{proof}
\end{lemma}

\begin{corollary}\label{corollary_Vshapedinternalzero}
If zero vertex weights are allowed in Lemmas \ref{lemma_Vshapedpendent} and
\ref{lemma_Vshapedinternal}, then optimal caterpillar $T\in
\mathcal{WT}(\mathbf{w}, \mathbf{d})$ with backbone $v_1, ...,v_q$ exists, where the following sequences are V-shaped for
$k=1,...,q$: internal vertex weights $\mu_T(v_k)$, internal vertex degrees $d_T(v_k)$, and pendent vertex weights $\mu_T(u_k)$ for any $u_k\in A_T(k)$.
\end{corollary}
The proof is similar to that of Corollary \ref{corollary_caterpillarzero}.

\section{Upper-Bound Estimate}\label{section_upper_bound}

Let us consider non-negative weight sequence $\mathbf{w}=(\mu_1, ..., \mu_n)$ and
natural degree sequence $\mathbf{d}=(d_1, ..., d_n)$, $d_i>1$, $i=1,...,q$, $d_i=1$,
$i=q+1,...,n$, $q\ge 4$. As before, sequence $\mathbf{d}$ is non-increasing, and
elements in $\mathbf{w}$ that share the same degree in $\mathbf{d}$, go in the
decreasing order of their weights. In this section it is shown that if weight sequence
$\mathbf{w}$ is monotone in degrees $\mathbf{d}$, then, for any $T\in
\mathcal{WT}(\mathbf{w},\mathbf{d})$ inequality $VWWI(T) \le UB(\mathbf{w},
\mathbf{d})$ holds, where
\begin{multline}\label{eq_UB}
UB(\mathbf{w},\mathbf{d}) =\mu\left[\frac{q+1}{4}\mu+\sum_{i=1}^q\mu_i\right]-
\sum_{i=q+1}^n\mu_i^2-\sum_{k=1}^{\lceil\frac{q}{2}\rceil}2M_k\left[kM_k+2\sum_{l=1}^{k-1}lM_l\right],
\end{multline}
$\mu:=\sum_{i=1}^n\mu_i$,
$M_k:=\frac{1}{2}\left(\mu_{2k-1}+\mu_{2k}+\sum_{i=1+D_{k-1}}^{D_k}\mu_{q+i}\right)$
for $k=1,...,\lfloor\frac{q}{2}\rfloor$,
$M_{\lfloor\frac{q}{2}\rfloor+1}:=\frac{1}{2}\left(\mu_q+\sum_{i=0}^{d_q-3}\mu_{n-i}\right)$,
 $D_0:=0$, $D_1:=d_1+d_2-2$,
$D_k=2+\sum_{i=1}^{2k}(d_i-2)$ for $k=2,...,\lceil\frac{q}{2}\rceil$.

According to expression (\ref{eq_VWWIass}), the optimal caterpillar problem (OCP)
 reduces to the assignment of internal and
pendent vertices to positions on the caterpillar's backbone taking into account
vertex degree constraints. Let $T\in \mathcal{WT}(\mathbf{w},\mathbf{d})$ be a
caterpillar with vertices $v_1, ..., v_n$ indexed so as $\mu_T(v_i)=\mu_i,
d_T(v_i)=d_i$, and let $X$ be the \emph{assignment} matrix, i.e., $x_{ik}=1$ when
vertex $v_i$, $i=1,...,n$, is associated with backbone position $k=1,...,q$, i.e.,
$v_i\in A_T(k)$, otherwise $x_{ik}=0$. Then the weight associated with backbone
position $k=1,...,q$ of caterpillar $T$ is written as
\begin{equation}\label{eq_weightX}
w_k(X | \mathbf{w}, \mathbf{d}) = \sum_{i=1}^n \mu_ix_{ik}, k=1,...,q,
\end{equation}
and position prices are written as
\begin{equation}\label{eq_priceX}
p_k(X | \mathbf{w}, \mathbf{d})=\sum_{l=1}^q w_l(X| \mathbf{w},
\mathbf{d})|k-l|=\sum_{i=1}^n \mu_i\sum_{l=1}^q x_{il} |k-l|, k=1,...,q.
\end{equation}
Taking into account expressions (\ref{eq_VWWIass}) and (\ref{eq_priceX}), 
the Wiener index of tree $T$ is written using matrix $X$ as
\begin{equation}\label{eq_VWWIX}
VWWI(X | \mathbf{w})=\frac{1}{2}\sum_{i,j=1}^n\sum_{k,l=1}^q\mu_i x_{ik}\mu_j
x_{jl}|k-l|+\sum_{i=1}^n\mu_i\sum_{j=1}^q\mu_j-\sum_{i=q+1}^n\mu_i^2,
\end{equation}
To simplify notation let us skip $\mathbf{w}$ and $\mathbf{d}$ when it does not
lead to confusion.

Only the first term depends on assignment matrix $X$, so OCP reduces to the
following binary quadratic program:

\begin{equation}
\begin{aligned}\label{eq_QAP}
\text{Maximize } & \frac{1}{2}\sum_{i,j=1}^n\sum_{k,l=1}^q\mu_i x_{ik}\mu_j
x_{jl}|k-l|
\end{aligned}
\end{equation}
subject to constraints:
\begin{align}
\text{discreteness: }& x_{ik} \in \{0,1\},i=1,...,n, k=1,...,q,\label{eq_zeroone}\\
\text{unique assignment: }& \sum_{k=1}^q x_{ik} = 1, i=1,...,n,\label{eq_assign}\\
\text{internal vertex assignment: }& \sum_{i=1}^q x_{ik} = 1, k=1,...,q,\label{eq_internal}\\
\text{balance of vertex degrees: }& \sum_{i=1}^n (d_i-2)x_{ik} = 0, k=2,...,q-1, \\
                                  & \sum_{i=1}^n (d_i-2)x_{ik} = -1, k=1,q.\label{eq_cluster}
\end{align}

This problem is similar to the Wiener-type QAP introduced in \cite{ccela2011wiener},
and is NP-complete, since it is a generalization of the classic partition problem
(Setting $q=2$, $n=2+2k$, $\mu_1=\mu_2=0$, $d_1=d_2=k$ makes a partition problem for
$2k$ elements.)

The continuous relaxation of problem (\ref{eq_QAP}) is obtained by replacing
discreteness constraints (\ref{eq_zeroone}) with box constraints
\begin{equation}\label{eq_unitrange}
x_{ik}\in[0,1].
\end{equation}

The function in expression (\ref{eq_QAP}) is concave on the feasible set: to verify,
use condition (\ref{eq_assign}) and obtain an equivalent problem
$$\frac{1}{2}\sum_{i,j=1}^n\sum_{k,l=1}^q\mu_i x_{ik}\mu_j
x_{jl}|k-l|=\frac{q-1}{2}\left(\mu_1+...+\mu_n\right)^2 - \frac{1}{2}\min_X
\mathbf{w}^\textrm{T}XPX^\textrm{T}\mathbf{w},$$ where matrix $P :=
(q-1-|k-l|)_{k,l=1}^q$ is positive semidefinite, since it can be written as the
following sum of positively definite rank-one matrices
\begin{multline}P=\sum_{i=1}^{q-1}\left[(\underbrace{1,...,1}_{i\text{ times}}, \underbrace{0,...,0}_{q-i\text{ times}})^\mathrm{T}\cdot(\underbrace{1,...,1}_{i\text{ times}}, \underbrace{0,...,0}_{q-i\text{ times}})+\right.\\
\left.+(\underbrace{0,...,0}_{i\text{ times}}, \underbrace{1,...,1}_{q-i\text{
times}})^\mathrm{T}\cdot(\underbrace{0,...,0}_{i\text{ times}},
\underbrace{1,...,1}_{q-i\text{ times}})\right].\nonumber
\end{multline}

Therefore, the relaxed OCP (ROCP) (\ref{eq_QAP}) with constraints
(\ref{eq_assign})-(\ref{eq_unitrange}) is a convex quadratic program with linear
constraints.


It is clear from (\ref{eq_QAP}) that OCP and ROCP do not change if one reverses the
numbering of caterpillar backbone positions. Therefore, if matrix $X$ is an optimal
solution of ROCP, then matrix $X'$ obtained from $X$ by applying inverse column
order is also optimal. Matrix $X''=\frac{X+X'}{2}$ is a feasible ROCP solution and,
since the optimization criterion is concave, $X''$ is also optimal. Therefore, when
studying ROCP, we can limit ourselves to \emph{symmetric solutions}, where
$x_{ik}=x_{i,q-k+1}$. Prices $p_k(X)$ are also symmetric with respect to the central
backbone position.

Similar to Lemma \ref{lemma_Vshapedpendent}, in the optimal ROCP solution vertices
with bigger weights are associated with outermost backbone positions. Lemma
\ref{lemma_pendentX} concerns pendent vertices, while Lemma \ref{lemma_internalX} is
about internal vertices.

\begin{lemma}\label{lemma_pendentX}
Let $\mathbf{w}$ be a weight sequence with \textbf{pairwise distinct} elements, and
let $X$ be an optimal symmetric solution of ROCP for weight sequence $\mathbf{w}$
and degree sequence $\mathbf{d}$. If $x_{i,k-1}>0$ and $x_{jk}>0$ for some
$i,j\in\{q+1,...,n\}$ and $k \le \lceil \frac{q}{2}\rceil$, then $i\le j$.
\begin{proof}
The proof is similar to that of Lemma \ref{lemma_Vshapedpendent}. Let us assume, by
contradiction, that $i > j$ (and hence, $\mu_j>\mu_i$, as weights are pairwise
distinct). Let us consider $n\times q$-matrix $\Delta$ with
$\Delta_{j,k-1}=\Delta_{ik}=-\Delta_{jk}=-\Delta_{i,k-1}=\delta>0$ and all other
elements being zeros. Assume that $\delta < \min\left(x_{i,k-1},x_{jk}\right)$, so
$X+\Delta$ is a feasible ROCP solution. Since $X$ is optimal,
\begin{equation}\label{ineq_VWWI_X}
VWWI(X+\Delta) - VWWI(X)= \delta\left(\mu_j-\mu_i\right)\left(p_{k-1}(X) - p_k(X) -
\delta(\mu_j-\mu_i)\right)\le 0.
\end{equation}

By assumption, $\mu_j-\mu_i>0$ and $\delta>0$. Hence, $p_{k-1}(X) - p_k(X) -
\delta(\mu_j-\mu_i)\le 0$, which is equivalent to $R - L + w_k(X) -
\delta(\mu_j-\mu_i)\le 0$, where $L:=\sum_{l=1}^{k-1} w_l(X)$, $R:=\sum_{l=k+1}^q
w_l(X)$. Since $X$ is a symmetric solution and $k \le \lceil \frac{q}{2}\rceil$,
inequality $R\ge L$ always holds and, consequently, $w_k(X) - \delta(w_j-w_i)\le 0$.
On the other hand, from $\delta < x_{jk}$ it follows that $w_k(X)\ge x_{jk}\mu_j >
\mu_j\delta$, so we obtain $\delta\mu_i< 0$, which is impossible. The obtained
contradiction proves the lemma.
\end{proof}
\end{lemma}

\begin{lemma}\label{lemma_internalX}
Assume weight sequence $\mathbf{w}$, consisting of \textbf{pairwise distinct}
elements, is monotone with respect to degree sequence $\mathbf{d}$. Let $X$ be an
optimal symmetric ROCP solution for weight sequence $\mathbf{w}$ and degree sequence
$\mathbf{d}$. If $x_{i,k-1}>0$ and $x_{jk}>0$ for some $i,j\in\{1,...,q\}$ and $k
\le \lceil \frac{q}{2}\rceil$, then $i\le j$.
\begin{proof}
The proof is similar to that of Lemma \ref{lemma_Vshapedinternal}. Assume, by
contradiction, that $j<i$. Since weights are pairwise distinct, it follows that
$\mu_j>\mu_i$, and since weights are degree-monotone, it follows that $d_j\ge d_i$.
Let us denote $d := d_j - d_i\ge 0$ and consider $n\times q$-matrix $\Delta$ with
$\Delta_{j,k-1}=\Delta_{ik}=-\Delta_{jk}=-\Delta_{i,k-1}=\delta>0$,
$\Delta_{v,k-1}=-\Delta_{vk}=\delta d x_{vk}(\sum_{v=1}^n x_{vk})^{-1},
v=q+1,...,n$, and all other elements being zeros; also introduce the shorthand
$\mu:=\delta d \sum_{v=q+1}^n \mu_v x_{vk}\left(\sum_{v=q+1}^n x_{vk}\right)^{-1}$.

Let us choose $\delta$ so as $\delta < \min\left(x_{i,k-1},x_{jk},\right)$ and
$d\delta < \sum_{v=q+1}^n x_{vk}$. Matrix $X+\Delta$ is a feasible ROCP solution, so
\begin{multline}VWWI(X+\Delta) - VWWI(X)=\\
 \delta\left(\mu+\mu_j-\mu_i\right)\left(p_{k-1}(X) - p_k(X) - \delta(\mu+\mu_j-\mu_i)\right).\end{multline}
The rest of the proof repeats that of Lemma \ref{lemma_pendentX}.
\end{proof}
\end{lemma}

Lemmas \ref{lemma_pendentX}
 and \ref{lemma_internalX} are refined by the following corollaries.

\begin{corollary}\label{corollary_exclusive_internal}
Under the conditions of Lemma \ref{lemma_internalX}, if $x_{ik}>0$ for some
$i\in\{1,...,q\}$ and $k\le \lceil \frac{q}{2}\rceil$, then $x_{jl}=0$ for all
$j=1,...,i$, $l=k+1,...,\lceil \frac{q}{2}\rceil$.
\begin{proof}
The proof is by contradiction. Firstly, let us assume that $x_{jl}>0$ for some
$j\in\{1,...,i-1\}$, $l\in\{k+1,...,\lceil \frac{q}{2}\rceil\}$.
From condition (\ref{eq_internal}), for every $k'\in \{k,...,l\}$ we have
$x_{i'k'}>0$ for some $i'\in\{1,...,q\}$. But by Lemma \ref{lemma_internalX}, if
$x_{i'k'}>0$ then $x_{i'',k'+1}=0$ for all $i''=1, ..., i'-1$, so $x_{jl}>0$ is
impossible for $1\le j<i$.

Secondly, let us assume that $x_{ik},x_{il}>0$ for some $i\in\{1,...,q\}$, $k,l\le
\lceil \frac{q}{2}\rceil$. Without loss of generality assume that $k<l$ and
$x_{jk}=0$ for all $j=1, ..., i-1$. From condition (\ref{eq_assign}) it follows that
$x_{ik}+x_{il}\le 1$, so $x_{ik}\le 1-x_{il}<1$. Then, from condition
(\ref{eq_cluster}), there exists such $j>i$ that $x_{jk}>0$. But it is shown above
that $x_{jk}>0, x_{il}>0$ is impossible. The obtained contradiction completes the
proof.
\end{proof}
\end{corollary}

\begin{corollary}\label{corollary_internalXzero}
If equal weights are allowed for weight sequence $\mathbf{w}$ in Lemma
\ref{lemma_internalX}, then such a symmetric ROCP solution $X$ \textbf{exists} that
$x_{ik}>0$ for some $i\in\{1,...,q\}$ and $k\in \{1,...,\lceil\frac{q}{2}\rceil\}$
implies that $x_{jl}=0$ for all $j=1,...,i$, $l=k+1,...,\lceil \frac{q}{2}\rceil$.
\begin{proof}
Let $\varepsilon:=\min_{i,j:\mu_i>\mu_j}\left(\mu_i-\mu_j\right)$ be the minimum
positive pairwise weight difference of the elements of sequence $\mathbf{w}$. Select
any positive $\alpha < \varepsilon$ and consider sequence
$\mathbf{w}^\alpha=(\mu^\alpha_i)_{i=1}^n=(\mu_i+\alpha/i)_{i=1}^n$ of pairwise
distinct positive weights. Let $X^\alpha$ be any optimal ROCP solution for weight
sequence $\mathbf{w}^\alpha$ and degree sequence $\mathbf{d}$. The feasible set of
ROCP does not depend on weights, so $X^\alpha$ is also feasible in ROCP for weight
sequence $\mathbf{w}$.

Let us consider an infinite sequence $X^\frac{\varepsilon}{t+1}$, $t=1,2, ...$. Its
elements belong to the bounded compact set, so, without loss of generality, it
converges (say, in Manhattan metric) to some feasible solution $X^*=(x_{ik}^*)$. It
is clear that $VWWI(X^\alpha|\mathbf{w})\le VWWI^*\le
VWWI(X^\alpha|\mathbf{w}^\alpha)\le VWWI(X^\alpha|\mathbf{w})+\alpha n^3$ for any
positive $\alpha$, where $VWWI^*$ is the optimal value of ROCP. Therefore,
$VWWI(X^*)=VWWI^*$ and $X^*$ is optimal.

Assume, by contradiction, that such $i,j\in \{1, ..., q\}$ and $k,l\le \lceil
\frac{q}{2}\rceil$, $j\le i$, $k<l$, exist that $x_{ik}^*>0$,$x_{jl}^*>0$, and
$i>j$. By Corollary \ref{corollary_exclusive_internal}, for any $t=1,2, ...$, and
$\alpha=\frac{\varepsilon}{t+1}$ either $x_{i,k}^\alpha=0$ or $x_{jl}^\alpha=0$, so
Manhattan distance between $X^\alpha$ and $X^*$ is at least $\min[x_{ik}^*,
x_{jl}^*]$ and sequence $X^\frac{\varepsilon}{t+1}$ cannot converge to $X^*$. The
obtained contradiction completes the proof.
\end{proof}
\end{corollary}

\begin{corollary}\label{corollary_exclusive_pendent}
Under the conditions of Lemma \ref{lemma_pendentX}, if $x_{ik}>0$ for some
$i\in\{q+1,...,n\}$ and $k\le \lceil \frac{q}{2}\rceil$, then $x_{jl}=0$ for all
$j=q+1,...,i$, $l=k+1,...,\lceil \frac{q}{2}\rceil$.
\end{corollary}
The proof is similar to that of Corollary \ref{corollary_exclusive_internal}.

\begin{corollary}\label{corollary_pendentXzero}
If equal weights are allowed in weight sequence $\mathbf{w}$ in Lemma
\ref{lemma_pendentX}, then such symmetric solution $X$ of ROCP exists that
$x_{ik}>0$ for some $i\in\{q+1,...,n\}$ and $k\in \{1,...,\frac{q}{2}\rceil\}$
implies that $x_{jl}=0$ for all $j=q+1,...,i$, $l=k+1,...,\lceil \frac{q}{2}\rceil$.
\end{corollary}
The proof is similar to that of Corollary \ref{corollary_internalXzero}.

\begin{theorem}
For odd $q$ the ``double V-shaped'' matrix
\begin{equation}\label{eq_ROCPX}\small
X = \begin{array}{cc} \begin{array}{ccccccc}
    &  &  & \lceil\frac{q}{2}\rceil &  &  &  \\
\end{array} &\\
\left(%
\begin{array}{ccccccccc}
  0.5 & 0 & ... & & & & ... & 0 & 0.5 \\
  0.5 & 0 & ... & & & & ... & 0 & 0.5 \\
  0 & 0.5 & ... & & & & ... & 0.5 & 0 \\
  0 & 0.5 & ... & & & & ... & 0.5 & 0 \\
  ... & ... & ... & ... & ... & ... & ... & ... & ... \\
  0 & 0 & ... & 0.5 &  0  & 0.5 & ... & 0 & 0 \\
  0 & 0 & ... & 0.5 &  0  & 0.5 & ... & 0 & 0 \\
  0 & 0 & ... & 0 & 1 & 0 & ... & 0 & 0 \\
  0.5 & 0 & ... & & & & ... & 0 & 0.5 \\
  ... & ... & ... & & & & ... & ... & .. \\
  0.5 & 0 & 0 & & & & ... & 0 & 0.5 \\
  0 & 0.5 & 0 & ...& & ... & 0 & 0.5 & 0 \\
  ... & ... & ... & & & & ... & ... & ... \\
  0 & 0.5 & 0 & ... & & ... & 0 & 0.5 & 0 \\
  ... & ... & ... & ... & ... & ... & ... & ... & ... \\
  0 & ... & 0 & 0.5 & 0 & 0.5 & 0 & ... & 0 \\
  ... & ... & ... & ... & ... & ... & ...& ... & ... \\
  0 & ... & 0 & 0.5 & 0 & 0.5 & 0 & ... & 0 \\
  0 & 0 & ... & 0 & 1 & 0 & ... & 0 & 0 \\
  ... & ... & ... & ... & ... & ... & & ... & ... \\
  0 & ... & ... & 0 & 1 & 0 & ... & ... & 0
\end{array}
\right) & \begin{array}{l} \left.
\begin{array}{c}
    \\
    \\
\end{array}
\right\}2\text{ times}\\
\left.
\begin{array}{c}
    \\
    \\
\end{array}
\right\}2\text{ times}\\
\\
\left.
\begin{array}{c}
    \\
    \\
\end{array}
\right\}2\text{ times}\\
\\
\left.
\begin{array}{c}
    \\
    \\
    \\
\end{array}
\right\}d_1+d_2-2\text{ times}\\
\left.
\begin{array}{c}
    \\
    \\
    \\
\end{array}
\right\}d_3+d_4-4\text{ times}\\
\\
\left.
\begin{array}{c}
    \\
    \\
    \\
\end{array}
\right\}d_{q-2}+d_{q-3}-4\text{ times}\\
\left.
\begin{array}{c}
    \\
    \\
    \\
\end{array}
\right\}d_q-2\text{ times}
\end{array}
\end{array}\normalsize
\end{equation}
is a solution of ROCP. In case of even $q$ the ``central'' column \emph{(}the one
marked with $\lceil\frac{q}{2}\rceil$\emph{)} is missing in (\ref{eq_ROCPX}).
\begin{proof}
Let us consider the case of odd $q$. From Corollaries \ref{corollary_internalXzero}
and \ref{corollary_pendentXzero} it follows that a symmetric ROCP soluton $X$ exists
where the first $\lceil \frac{q}{2}\rceil$ columns have a single non-zero element in
each row. From the symmetry of solution $X$ and condition (\ref{eq_assign}), non-zero
element $x_{ik}=\frac{1}{2}$ for $k<\lceil \frac{q}{2}\rceil$ and $x_{i\lceil
\frac{q}{2}\rceil}=1$. Corollary \ref{corollary_internalXzero} says that non-zero
elements for the first $q$ rows and first $\lceil \frac{q}{2}\rceil$ columns go from
the top-left corner to the bottom-right, and from condition (\ref{eq_internal}) it
follows that column $\lceil \frac{q}{2}\rceil$ has a single non-zero element, while
all the other columns have exactly two non-zero elements in the first $q$ rows.

The placement of non-zero elements for rows from $q+1$ to $n$ is justified in the
same manner, with the only difference that the number of non-zero elements in each
column is determined by condition (\ref{eq_cluster}). The case of even $q$ is
considered similarly, except that column $\lceil \frac{q}{2}\rceil$ does not play
any special role.
\end{proof}
\end{theorem}

\begin{corollary}
If weight sequence $\mathbf{w}$ is monotone in degree sequence $\mathbf{d}$,
expression (\ref{eq_UB}) gives an upper bound of the Wiener index over
$\mathcal{WT}(\mathbf{w},\mathbf{d})$.
\begin{proof}
Expression (\ref{eq_UB}) is obtained from expressions (\ref{eq_VWWIX}) and
(\ref{eq_ROCPX}) with a series of simplifications that employ symmetry of function
(\ref{eq_QAP}) and of solution $X$. Expression (\ref{eq_VWWIX}) gives an upper bound
of the Wiener index for optimal caterpillar since matrix $X$ is a solution of a
continuous relaxation of OCP.
\end{proof}
\end{corollary}

\section{Quality of Upper Bound}\label{section_quality}

Let us show first that upper bound (\ref{eq_UB}) is tight.
\begin{theorem}\label{th_tight}
If in weight sequence $\mathbf{w}$ being monotone in degree sequence $\mathbf{d}$ we
have $\mu_{2i-1}=\mu_{2i}$, $d_{2i-1}=d_{2i}$ for $i=1,..., \lfloor
\frac{q}{2}\rfloor$ and $\mu_{q+2i-1}=\mu_{q+2i}$ for $i=1,..., \lfloor
\frac{n-q}{2}\rfloor$ \emph{(}i.e., the number of vertices with distinct degree and
weight is even, while the pendent/internal vertex with the smallest weight may be
unpaired in case of the odd number of pendent/internal vertices\emph{)}, then
$UB(\mathbf{w},\mathbf{d})=VWWI(T)$ for some $T\in
\mathcal{WT}(\mathbf{w},\mathbf{d})$.
\begin{proof}
Let $X=(x_{ik})$ be a ROCP solution defined by expression (\ref{eq_ROCPX}) and
consider $X'=(x_{ik}')$, where for all $i=1,...,n$ and $k=1,...,q$
$$x_{ik}'=\left\{%
\begin{array}{ll}
    1, & \hbox{$x_{ik}>0,i\le q, k\le\lceil{q/2}\rceil\text{ and } i\text{ is odd}$;} \\
    1, & \hbox{$x_{ik}>0,i\le q, k\ge\lceil{q/2}\rceil\text{ and } i\text{ is even}$;} \\
    1, & \hbox{$x_{ik}>0,i> q, k\le\lceil{q/2}\rceil\text{ and } i-q\text{ is odd}$;} \\
    1, & \hbox{$x_{ik}>0,i> q, k\ge\lceil{q/2}\rceil\text{ and } i-q\text{ is even}$;} \\
    0, & \hbox{otherwise,} \\
\end{array}\right.$$ i.e., in every pair of vertices with the same weight and degree
distributed in $X$ equally between positions $k$ and $q-k+1$ of the caterpillar
backbone, one is assigned to position $k$, and the other is assigned to position
$q-k+1$. It is clear that $X'$ is a feasible OCP solution, and
$UB(\mathbf{w},\mathbf{d})=VWWI(X)=VWWI(X')$, which completes the proof.
\end{proof}
\end{theorem}
From the proof of Theorem \ref{th_tight} it follows that for odd $q$ its conditions
could be a bit weakened: $d_q-2$ pendent vertices with the smallest weight may be
unpaired.

Let $T^*(\mathbf{w},\mathbf{d})$ be an optimal caterpillar for weight sequence
$\mathbf{w}$ and degree sequence $\mathbf{d}$. The relative error of upper bound
(\ref{eq_UB}) is defined as
\begin{equation}\label{eq_RE}
RE(\mathbf{w},\mathbf{d}) :=
\frac{UB(\mathbf{w},\mathbf{d})}{VWWI(T^*(\mathbf{w},\mathbf{d}))}-1.
\end{equation}

Unfortunately, $RE(\mathbf{w},\mathbf{d})$ can be arbitrary high. For example, if
$\mu_1>0$ and $\mu_i=0$ for $i>1$, then $VWWI(T)=0$ for every $T\in
\mathcal{WT}(\mathbf{w},\mathbf{d})$, but $UB(\mathbf{w},\mathbf{d})>0$, so
$\frac{UB(\mathbf{w},\mathbf{d})}{VWWI(T^*(\mathbf{w},\mathbf{d}))}$ is unbounded.
Nevertheless, it is shown below that the average relative error is low enough,
especially for large trees. Let us define
$$\delta_k=\left\{%
\begin{array}{ll}
    $1$, & \hbox{$k\in \{1,q\}$;} \\
    $2$, & \hbox{otherwise.} \\
\end{array}\right.$$
The greedy heuristics shown in Listing \ref{lst:greedy} constructs a (nearly
optimal) caterpillar for weight sequence $\mathbf{w}$ and degree sequence
$\mathbf{d}$ by sequentially assigning vertices to the caterpillar backbone position
with the maximum current price. Prices are re-calculated after each iteration, so
two consequent vertices are typically spread as far as possible.

\begin{algorithm}[!ht]
\floatname{algorithm}{Listing} \caption{Greedy heuristic algorithm builds an
approximately optimal caterpillar for weight sequence $\mathbf{w}$ and degree
sequence $\mathbf{d}$}\label{lst:greedy}
\begin{algorithmic}[1]
\Function{GreedyCaterpillar}{$\mathbf{w}$, $\mathbf{d}$}  \State $X \gets 0_{n\times
q}$ \Comment{Start from the empty assignment table} \State $X(q+1,1) \gets 1$
\Comment{Assign $1^{st}$ pendent vertex to the outermost backbone position} \State
$i \gets 1$ \Comment{Internal vertex counter} \State $j = q+2$ \Comment{Pendent
vertex counter}

 \While {$i<q$ or $j<n$} \Comment{Loop until unassigned vertices exist}
\State
    $k \gets \argmax\{p_{k'}(X)| k':\sum_{i'=1}^q x_{i'k'}<1\}$\State
    $l \gets \argmax\{p_{l'}(X)| l':\sum_{j'=q+1}^n x_{j'l'}<\sum_{i'=1}^q
    d_{i'}x_{i'l'}-\delta_{l'}\}$

    \If {$\mu_ip_k(X)>\mu_jp_l(X)$} \State $x_{ik}\gets 1$ , $i\gets i+1$\Else
\State $x_{jl}\gets 1$, $j\gets j+1$\EndIf \EndWhile

\State \Return{$X$} \EndFunction
\end{algorithmic}
\end{algorithm}

\begin{figure}\center
  \includegraphics[width=0.8\textwidth]{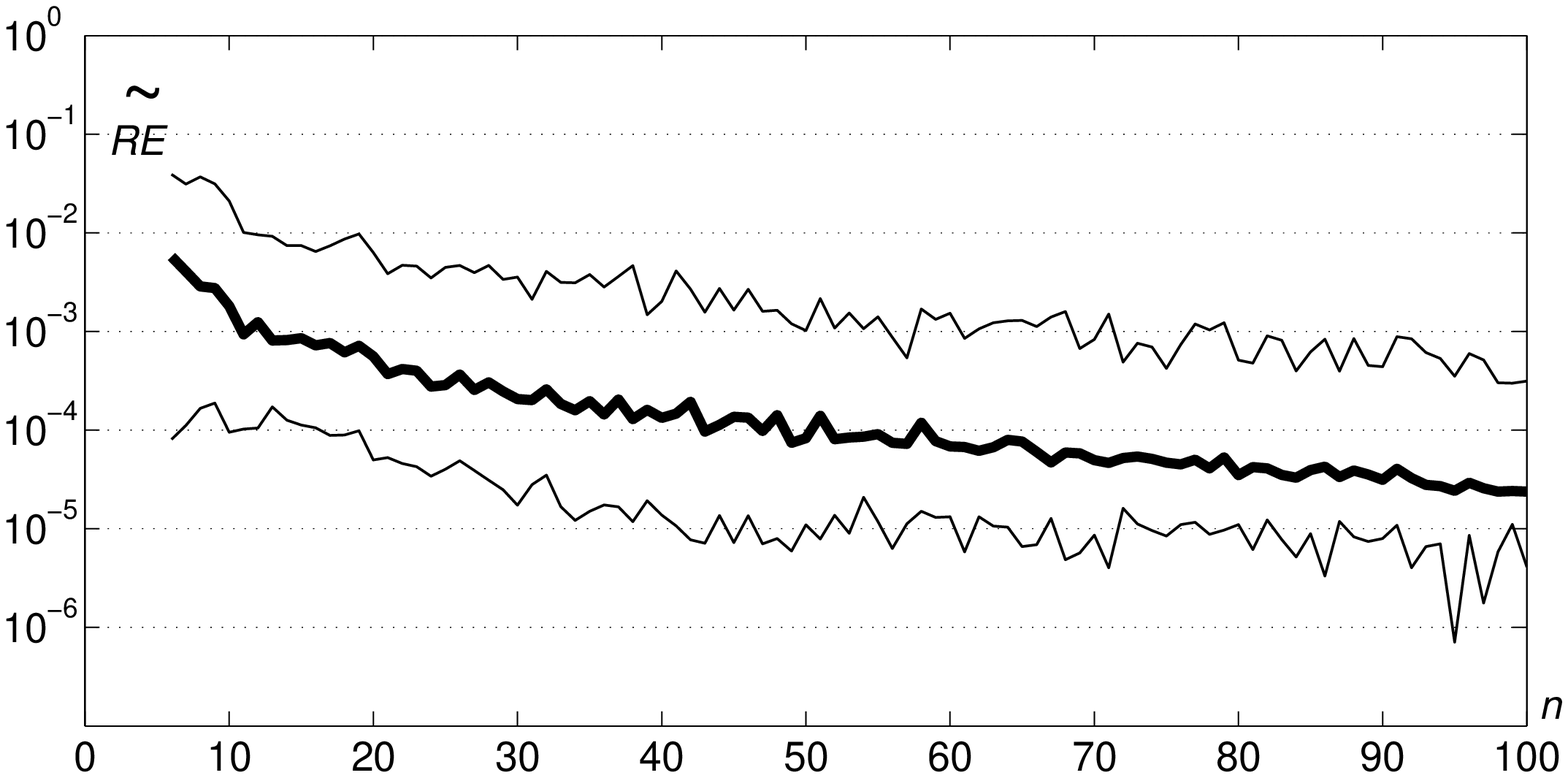}
\caption{Estimated relative error of upper bound (\ref{eq_UB}) \emph{vs} graph order
(logarithmic scale for the vertical axis)} \label{fig:UBError}
\end{figure}

The relative error can be estimated from above as $$RE(\mathbf{w},\mathbf{d}) \le
\widetilde{RE}(\mathbf{w},\mathbf{d})=
\frac{UB(\mathbf{w},\mathbf{d})}{VWWI(\text{\textproc{GreedyCaterpillar}}(\mathbf{w},\mathbf{d})|\mathbf{w})}-1.$$
In Figure \ref{fig:UBError} the value of $\widetilde{RE}(\cdot)$ is presented as a
function of the graph order. The minimum error is zero and the maximum tends to
infinity, but these are extremely rare events. Therefore, the median (the bold
line), the upper and the lower decile values (two thin lines) are shown in Figure
\ref{fig:UBError} for 1000 randomly generated weight and degree sequences for each
vertex count $n=6,...,100$. The median error never exceeds $1\%$, and even for
relatively small trees (those with $n>12$ vertices) nine of ten graphs have error
less than $1\%$. For moderately sized trees ($n>50$) the median relative error is
less than $0.01\%$, and for $90\%$ of trees the relative error does not exceed
$0.2\%$.

\section{Branch and Bound Scheme}\label{section_BandB}

Although it is shown in the previous section that both upper bound (\ref{eq_UB}) and
heuristic algorithm \textproc{GreedyCaterpillar} have good average quality, the
question of the exact solution of OCP is still open. Since OCP reduces to the convex
binary quadratic minimization program (\ref{eq_QAP}) with linear constraints
(\ref{eq_zeroone})-(\ref{eq_cluster}), the branch and bound algorithms implemented
in commercial solvers (like CPLEX) can be used to find an optimal caterpillar.
Nevertheless, due to the large search space ($2^{n\times q}$ elements), they show
low performance being inapplicable even for trees with a dozen of vertices.

In case of weights being monotone in degrees the search space can be decreased
dramatically taking into account the characterization of the structure of an optimal
caterpillar presented in Section \ref{section_structure}. Let us restrict attention
to optimal caterpillars with V-shaped weight sequences as stated by Corollary
\ref{corollary_Vshapedinternalzero}. Then it is clear, that one obtains any V-shaped
assignment of internal vertices to positions on the caterpillar backbone by
sequentially running through internal vertices in the order of descending weights
and deciding whether to assign the vertex to the rightmost or to the leftmost vacant
backbone position. Each decision is binary (right or left), so $2^q$ alternative
assignments are to be considered.

In the same manner, with positions of internal vertices being fixed, any V-shaped
assignment of pendent vertices is obtained by sequentially assigning pendent
vertices to the rightmost or to the leftmost vacant backbone position taking care
about vertex degree constraints. Therefore, the total search space has $2^n$
alternative assignments, which makes great economy.

The branch and bound algorithm presented in Listing \ref{lst:BandB} can be used to
solve OCP for trees of moderate order ($n \approx 30$) on a PC. It systematically
explores the search space by recursively assigning internal and pendent vertices to
backbone positions as explained above (the branching takes place on whether the
vertex is assigned to the left or to the right). If the best known solution
outperforms the upper bound of $VWWI$ given the current partially fixed solution,
then further branching is not required (enumeration is bounded).

\textproc{GreedyCaterpillar($\mathbf{w},\mathbf{d}$)} is used as a starting best
known solution. The upper bound $UB(\mathbf{w},\mathbf{d}|X_0)$ given partial
assignment matrix $X^0=(x_{ik}^0)$ is calculated as a solution of ROCP
(\ref{eq_QAP}), (\ref{eq_assign})-(\ref{eq_unitrange}) with additional box
constraints $x_{ik}=1$ for such $i$ and $k$ that $x_{ik}^0=1$. This convex quadratic
program is solved efficiently by many open-source and commercial packages. To
improve performance, in our Matlab implementation\footnote{The code is available
online at \texttt{http://www.mtas.ru/upload/maxWiener.zip}} branching on the next
pendent vertex occurs immediately as a pair of vacant left and right positions
appears on the caterpillar backbone, and the order of recursion (right then left or
vice versa) is driven by the current price $p_k(X_0)$ of backbone position, as in
Listing \ref{lst:greedy}.

\begin{algorithm}[!ht]
\floatname{algorithm}{Listing} \caption{Branch and bound algorithm for
OCP}\label{lst:BandB}
\begin{algorithmic}[1]
\Function{OptimalCaterpillar}{$\mathbf{w}$, $\mathbf{d}$} \State $X_0 \gets
0_{n\times q}$ \Comment{Start from empty assignment table} \State $Best \gets $
GreedyCaterpillar($\mathbf{w}, \mathbf{d}$) \Comment{The best known solution} \State
Branching($X_0$) \Comment{Recursive solution branching} \State Return $Best$
\Comment{Branching leaves the solution in $Best$} \EndFunction
\Function{Branching}{$X_0$} \Comment{Recursive branching starting from partial
solution $X_0$} \State $i \gets \min\{j:\sum_{l=1}^q x_{jl}^0=0\}$ \Comment{$1^{st}$
unassigned vertex in $X_0$} \If {$i<q$} \Comment{If unassigned internal vertex
exists} \State $k \gets \max\{l:\sum_{j=1}^q x_{jl}^0=0\}$, $X_R\gets X_0$,
$x_{ik}^R \gets 1$ \Comment{Assign internal right} \If
{$VWWI$($UB(\mathbf{w},\mathbf{d}|X_R)$)$>VWWI(Best)$} \State Branching($X_R$)
\EndIf \State $k \gets \min\{l:\sum_{j=1}^q x_{jl}^0=0\}$, $X_L\gets X_0$, $x_{ik}^L
\gets 1$ \Comment{Assign internal left} \If
{$VWWI$($UB(\mathbf{w},\mathbf{d}|X_L)$)$>VWWI(Best)$} \State Branching($X_L$)
\EndIf
 \ElsIf {$i \le n$} \Comment{If unassigned pendent vertices exist}
 \State $k \gets \max\{l:\sum_{j=q+1}^n x_{jl}^0<\sum_{j=1}^q
    d_jx_{jl}^0-\delta_l\}$, $X_R\gets X_0$,
$x_{ik}^R \gets 1$ \Comment{Assign right} \If
{$VWWI$($UB(\mathbf{w},\mathbf{d}|X_R)$)$>VWWI(Best)$} \State Branching($X_R$)
\EndIf \State $k \gets \min\{l:\sum_{j=q+1}^n x_{jl}^0<\sum_{j=1}^q
    d_jx_{jl}^0-\delta_l\}$, $X_L\gets X_0$,
$x_{ik}^L \gets 1$ \Comment{Assign left} \If
{$VWWI$($UB(\mathbf{w},\mathbf{d}|X_L)$)$>VWWI(Best)$} \State Branching($X_L$)
\EndIf \ElsIf {$VWWI(X_0)
> VWWI(Best)$} $Best \gets X_0$ \Comment{Update the best}
 \EndIf
\EndFunction
\end{algorithmic}
\end{algorithm}

\begin{figure}\center
  \includegraphics[width=0.8\textwidth]{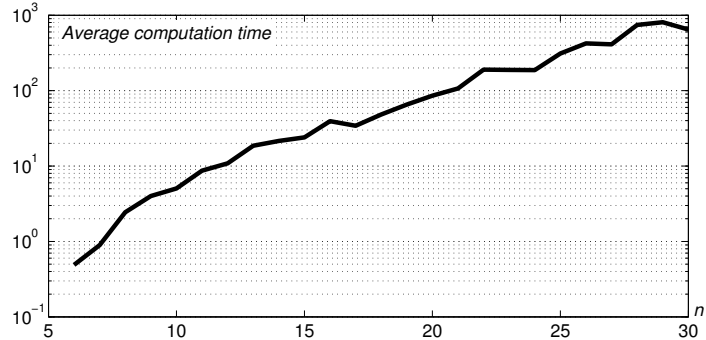}
\caption{Average computation time of the branch and bound algorithm \emph{vs} graph
order (logarithmic scale for the vertical axis)} \label{fig:BBtime}
\end{figure}

To evaluate the performance of the algorithm, 100 random weight and degree sequences
were generated for every graph order $n = 6, ..., 30$. In Figure \ref{fig:BBtime}
the average computation time is shown as a function of $n$. It grows exponentially,
which is expectable for the NP-complete problem.

\begin{figure}\center
  \includegraphics[width=0.8\textwidth]{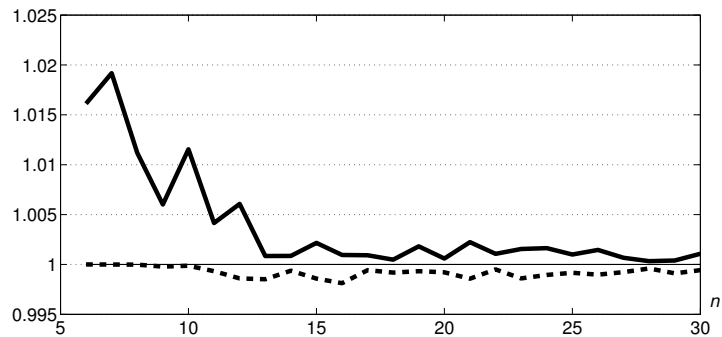}
\caption{Upper bound (\ref{eq_UB}) (thick solid line) and the value of the Wiener
index for \textproc{GreedyCaterpillar} (dashed line) with respect to the Wiener
index for the optimal caterpillar (thin horizontal line) \emph{vs} graph order $n$}
\label{fig:GC_BB_UB}
\end{figure}

In Figure \ref{fig:GC_BB_UB} upper bound (\ref{eq_UB}) and the value of the Wiener
index for \textproc{GreedyCaterpillar} are compared with the optimal Wiener index
value. The thick solid line in Figure \ref{fig:GC_BB_UB} is
$RE(\mathbf{w},\mathbf{d})+1=\frac{UB(\mathbf{w},\mathbf{d})}{VWWI(T^*(\mathbf{w},\mathbf{d}))}$
averaged over 100 random weight and degree sequences for the given graph order $n$,
while the dashed line is
$\frac{VWWI(\text{\textproc{GreedyCaterpillar}}(\mathbf{w},\mathbf{d}))}{VWWI(T^*(\mathbf{w},\mathbf{d}))}$,
also averaged. It can be seen that, although the gap between the index value for the greedy
tree and the upper bound is typically small, filling this small gap by the branch
and bound algorithm may take much time.

\section{Conclusion}

This article contributes to the literature on the Wiener index by studying the
Wiener index maximization problem over the set $\mathcal{WT}(\mathbf{w},\mathbf{d})$
of trees with the given vertex weight sequence $\mathbf{w}$ and degree sequence
$\mathbf{d}$. The results of \cite{schmuck2012greedy} were extended to the Wiener
index for graphs with weighted vertices: it was proven that if vertex weight
sequence $\mathbf{w}$ is monotone in degree sequence $\mathbf{d}$, then there is an
optimal caterpillar with internal/pendent vertex weights monotonously increasing
from some central point to the ends of its backbone.

For weight sequences being monotone in degrees, closed-form expression (\ref{eq_UB})
was proposed for the upper bound of the Wiener index value over
$\mathcal{WT}(\mathbf{w},\mathbf{d})$. It was shown that this upper bound is tight,
and an efficient greedy heuristics was proposed that approximates well the optimal
tree. Finally, a branch and bound scheme was proposed for the exact solution of this
NP-complete problem and computational analysis of its performance was accomplished.

Most of the results of this article are limited to the case of weight sequences being
monotone in degrees, when the weight of an internal vertex does not decrease in its
degree (no restrictions are imposed on weights of pendent vertices). The general
case of non-monotonous weight sequences seems more complicated. Corollary
\ref{corollary_caterpillarzero} says that an optimal caterpillar exists, but weight
sequences do not have to be V-shaped, so, expression (\ref{eq_UB}) for the upper
bound is directly inapplicable, although the solution of relaxed OCP (\ref{eq_QAP}),
(\ref{eq_assign})-(\ref{eq_unitrange}) still gives an efficiently calculable
upper-bound estimate, and the greedy algorithm still builds some suboptimal trees.

  \bibliography{GoubkoWiener}

\begin{thebibliography}{10}
\expandafter\ifx\csname url\endcsname\relax
  \def\url#1{\texttt{#1}}\fi
\expandafter\ifx\csname urlprefix\endcsname\relax\def\urlprefix{URL }\fi
\expandafter\ifx\csname href\endcsname\relax
  \def\href#1#2{#2} \def\path#1{#1}\fi

\bibitem{dobrynin2001wiener}
A.~A. Dobrynin, R.~Entringer, I.~Gutman, {W}iener index of trees: theory and
  applications, Acta Appl. Math. 66~(3) (2001) 211--249.
\newblock \href {http://dx.doi.org/10.1023/A:1010767517079}
  {\path{doi:10.1023/A:1010767517079}}.

\bibitem{wang2008extremal}
H.~Wang, The extremal values of the {W}iener index of a tree with given degree
  sequence, Discrete App. Math. 156~(14) (2008) 2647--2654.

\bibitem{zhang2008wiener}
X.-D. Zhang, Q.-Y. Xiang, L.-Q. Xu, R.-Y. Pan, The {W}iener index of trees with
  given degree sequences, MATCH Commun. Math. Comput. Chem. 60~(2) (2008)
  623--644.

\bibitem{klavvzar1997wiener}
S.~Klav{\v{z}}ar, I.~Gutman, {W}iener number of vertex-weighted graphs and a
  chemical application, Discrete Appl. Math. 80~(1) (1997) 73--81.

\bibitem{kelenc2015edge}
A.~Kelenc, S.~Klav{\v{z}}ar, N.~Tratnik, The edge-{W}iener index of benzenoid
  systems in linear time, MATCH Commun. Math. Comput. Chem. 74~(3) (2015)
  521--532.

\bibitem{gao2015vertex}
W.~Gao, W.~Wang, The vertex version of weighted {W}iener number for bicyclic
  molecular structures, Comp. and Math. Meth. Med. 2015.
\newblock \href {http://dx.doi.org/10.1155/2015/418106}
  {\path{doi:10.1155/2015/418106}}.

\bibitem{goubkomiloserdov2016wiener}
M.~Goubko, O.~Miloserdov, Simple alcohols with the lowest normal boiling point
  using topological indices, MATCH Commun. Math. Comput. Chem. 75 (2016)
  29--56.

\bibitem{goubko2016wiener}
M.~Goubko, Minimizing {W}iener index for vertex-weighted trees with given
  weight and degree sequences, MATCH Commun. Math. Comput. Chem. 75 (2016)
  3--27.

\bibitem{shi1993average}
R.~Shi, The average distance of trees, Sys. Sci. Math. Sci. 6 (1993) 18--24.

\bibitem{ccela2011wiener}
E.~{\c{C}}ela, N.~S. Schmuck, S.~Wimer, G.~J. Woeginger, The {W}iener maximum
  quadratic assignment problem, Discrete Optim. 8~(3) (2011) 411--416.

\bibitem{wiener1947structural}
H.~Wiener, Structural determination of paraffin boiling points, J. Amer. Chem.
  Soc. 69~(1) (1947) 17--20.

\bibitem{knor2016mathematical}
M.~Knor, R.~{\v{S}}krekovski, A.~Tepeh, Mathematical aspects of {W}iener index,
  Ars Math. Contemp. 11~(2) (2016) 327--352.

\bibitem{diudea1998wiener}
M.~V. Diudea, I.~Gutman, {W}iener-type topological indices, Croatica Chem. Acta
  71~(1) (1998) 21--51.

\bibitem{gutman1993some}
I.~Gutman, Y.-N. Yeh, S.-L. Lee, Y.-L. Luo, Some recent results in the theory
  of the {W}iener number, Indian J. Chem. 32 (1993) 651--661.

\bibitem{rouvray1987modeling}
D.~H. Rouvray, The modeling of chemical phenomena using topological indices, J.
  Comput. Chem. 8~(4) (1987) 470--480.

\bibitem{lang2016novel}
R.~Lang, T.~Li, D.~Mo, Y.~Shi, A novel method for analyzing inverse problem of
  topological indices of graphs using competitive agglomeration, Appl. Math.
  Comp. 291 (2016) 115--121.
\newblock \href {http://dx.doi.org/10.1016/j.amc.2016.06.048}
  {\path{doi:10.1016/j.amc.2016.06.048}}.

\bibitem{fischermann2002wiener}
M.~Fischermann, A.~Hoffmann, D.~Rautenbach, L.~Sz{\'e}kely, L.~Volkmann,
  {W}iener index versus maximum degree in trees, Discrete Appl. Math. 122~(1)
  (2002) 127--137.

\bibitem{lin2016segment}
H.~Lin, M.~Song, On segment sequences and the {W}iener index of trees, MATCH
  Commun. Math. Comput. Chem. 75~(1) (2016) 81--89.

\bibitem{lin2014extremal}
H.~Lin, Extremal {W}iener index of trees with given number of vertices of even
  degree, MATCH Commun. Math. Comput. Chem. 72~(1) (2014) 311--320.

\bibitem{krnc2016wiener}
M.~Krnc, R.~\v{S}krekovski, On {W}iener inverse interval problem, MATCH Commun.
  Math. Comput. Chem. 75~(1) (2016) 71--80.

\bibitem{ramane2016note}
H.~S. Ramane, V.~V. Manjalapur, Note on the bounds on wiener number of a graph,
  MATCH Commun. Math. Comput. Chem. 76~(1) (2016) 19--22.

\bibitem{shi2016chemical}
L.~Shi, Chemical indices, mean distance, and radius, MATCH Commun. Math.
  Comput. Chem. 75~(1) (2016) 57--70.

\bibitem{merris1989edge}
R.~Merris, An edge version of the matrix-tree theorem and the {W}iener index,
  Linear Multilinear Algebra 25~(4) (1989) 291--296.

\bibitem{merris1990distance}
R.~Merris, The distance spectrum of a tree, J. Graph Theory 14~(3) (1990)
  365--369.

\bibitem{mohar1991eigenvalues}
B.~Mohar, Eigenvalues, diameter, and mean distance in graphs, Graphs and
  Combin. 7~(1) (1991) 53--64.

\bibitem{mohar1991laplacian}
B.~Mohar, The {L}aplacian spectrum of graphs, in: Y.~Alavi, G.~Chartrand,
  O.~Ollermann, A.~J. Schwenk (Eds.), Graph Theory, Combinatorics, and
  Applications, Vol.~2, Wiley, 1991, pp. 871--898.

\bibitem{indulal2009sharp}
G.~Indulal, Sharp bounds on the distance spectral radius and the distance
  energy of graphs, Linear Algebra Appl. 430~(1) (2009) 106--113.

\bibitem{gallai1960graphs}
T.~Gallai, P.~Erd\H{o}s, Graphs with prescribed degree of vertices
  ({H}ungarian), Mat. Lapok 11 (1960) 264--274.

\bibitem{burr1988extremal}
S.~A. Burr, P.~Erd\H{o}s, R.~J. Faudree, A.~Gy{\'a}rf{\'a}s, R.~Schelp,
  Extremal problems for degree sequences, Combin. 52 (1988) 183--193.

\bibitem{biyikoglu2008graphs}
T.~B{\i}y{\i}ko\u{g}lu, J.~Leydold, Graphs with given degree sequence and
  maximal spectral radius, Electron. J. Combin. 15 (2008) R119.

\bibitem{gutman2013degree}
I.~Gutman, Degree-based topological indices, Croatica Chem. Acta 86~(4) (2013)
  351--361.

\bibitem{gutman2004first}
I.~Gutman, K.~C. Das, The first {Z}agreb index 30 years after, MATCH Commun.
  Math. Comput. Chem. 50 (2004) 83--92.

\bibitem{gutman2015graphs}
I.~Gutman, M.~K. Jamil, N.~Akhter, Graphs with fixed number of pendent vertices
  and minimal first {Z}agreb index, Trans. Combin. 4~(1) (2015) 43--48.

\bibitem{das2015zagreb}
K.~C. Das, K.~Xu, J.~Nam, {Z}agreb indices of graphs, Front. Math. China 10~(3)
  (2015) 567--582.

\bibitem{ghorbani2012new}
M.~Ghorbani, M.~A. Hosseinzadeh, A new version of {Z}agreb indices, Filomat
  26~(1) (2012) 93--100.

\bibitem{eliasi2012multiplicative}
M.~Eliasi, A.~Iranmanesh, I.~Gutman, Multiplicative versions of first {Z}agreb
  index, MATCH Commun. Math. Comput. Chem. 68~(1) (2012) 217--230.

\bibitem{fath2011old}
G.~Fath-Tabar, Old and new {Z}agreb indices of graphs, MATCH Commun. Math.
  Comput. Chem. 65 (2011) 79--84.

\bibitem{xu2014maximizing}
K.~Xu, K.~C. Das, S.~Balachandran, Maximizing the {Z}agreb indices of $(n,
  m)$-graphs, MATCH Commun. Math. Comput. Chem. 72 (2014) 641--654.

\bibitem{furtula2015forgotten}
B.~Furtula, I.~Gutman, A forgotten topological index, J. Math. Chem. 53~(4)
  (2015) 1184--1190.

\bibitem{su2015graphs}
G.~Su, J.~Tu, K.~C. Das, Graphs with fixed number of pendent vertices and
  minimal {Z}eroth-order general {R}andi{\'c} index, Appl. Math. Comput. 270
  (2015) 705--710.

\bibitem{xu2012unified}
K.~Xu, H.~Hua, A unified approach to extremal multiplicative {Z}agreb indices
  for trees, unicyclic and bicyclic graphs, MATCH Commun. Math. Comput. Chem.
  68~(1) (2012) 241--256.

\bibitem{zhang2015extremal}
X.-D. Zhang, Extremal graph theory for degree sequences, arXiv preprint
  arXiv:1510.01903.

\bibitem{liu2015recent}
M.~Liu, B.~Liu, K.~C. Das, Recent results on the majorization theory of graph
  spectrum and topological index theory, Electron. J. Lin. Algebra 30~(1)
  (2015) 402--421.
\newblock \href {http://dx.doi.org/10.13001/1081-3810.3086}
  {\path{doi:10.13001/1081-3810.3086}}.

\bibitem{schmuck2012greedy}
N.~S. Schmuck, S.~G. Wagner, H.~Wang, Greedy trees, caterpillars, and
  {W}iener-type graph invariants, MATCH Commun. Math. Comput. Chem. 68~(1)
  (2012) 273--292.

\bibitem{gutman2009terminal}
I.~Gutman, B.~Furtula, M.~Petrovi{\'c}, Terminal {W}iener index, J. Math. Chem.
  46~(2) (2009) 522--531.

\bibitem{gutman2010survey}
I.~Gutman, B.~Furtula, A survey on terminal {W}iener index, novel molecular
  structure descriptors -- theory and applications i Edition, Univ. Kragujevac,
  Kragujevac, 2010, pp. 173--190.

\bibitem{mohar1988compute}
B.~Mohar, T.~Pisanski, How to compute the {W}iener index of a graph, J. Math.
  Chem. 2~(3) (1988) 267--277.

\bibitem{spielman2012}
D.~Spielman, Spectral graph theory, combin. sci. comp. Edition, CRC Press, Boca
  Raton, 2012, pp. 495--524.

\bibitem{huffman1952method}
D.~A. Huffman, A method for the construction of minimum-redundancy codes, Proc.
  IRE 40~(9) (1952) 1098--1101.

\end{thebibliography}

\end{document}